\documentclass[smallextended]{svjour3}
\usepackage{hyperref}
\smartqed
\usepackage{mathptmx,amsmath,amssymb,amsfonts,bm,mathtools,graphicx,color,cite}


\DeclareMathOperator*{\argmin}{arg\;min}

\DeclareMathOperator*{\argsup}{arg\;sup}

\newcommand{\rhorho}{\varphi}
\newcommand{\varphivarphi}{\rho}
\newcommand{\bx}{{\bf x}}

\newcommand{\vy}{{\bm y}}
\newcommand{\mmm}{{M}}
\newcommand{\kkk}{{K}}
\newcommand{\KKK}{{Z}}

\newcommand{\vertiii}[1]{{\left\vert\kern-0.25ex\left\vert\kern-0.25ex\left\vert #1 
    \right\vert\kern-0.25ex\right\vert\kern-0.25ex\right\vert}}
    
\journalname{Journal of Scientific Computing}

\begin{document}

\title{An improved discrete least-squares/reduced-basis\\ method for parameterized elliptic PDEs\thanks{This material is based upon work supported in part by the U.S. Air Force of Scientific Research under grants 1854-V521-12 and FA9550-15-1-0001; the U.S. Defense Advanced Research Projects Agency, Defense Sciences Office under contract and award HR0011619523 and 1868-A017-15; the U.S. Department of Energy, Office of Science, Office of Advanced Scientific Computing Research, Applied Mathematics program under contracts and awards ERKJ259, ERKJ320, DE-SC0009324, and DE-SC0010678; by the National Science Foundation under the award numbers 1620027 and 1620280.}
}

\titlerunning{An improved discrete least-squares/reduced-basis method for parameterized elliptic PDEs}

\author{Max Gunzburger \and Michael Schneier \\ Clayton Webster \and Guannan Zhang}

\authorrunning{M. Gunzburger, M. Schneier, C. Webster, and G. Zhang}

\institute{M. Gunzburger and M. Schneier\at
           Department of Scientific Computing, Florida State University, Tallahassee FL 32306-4120, USA.
           \email{mgunzburger@fsu.edu, mhs13c@my.fsu.edu} 
\and
           C. Webster and G. Zhang\at
           Computational and Applied Mathematics Division, Oak Ridge National Laboratory, Oak Ridge TN 37831-6164, USA.
           \email{webstercg@ornl.gov, zhangg@ornl.gov}              
}

\date{Received: date / Accepted: date}

\maketitle

\begin{abstract}
It is shown that the computational efficiency of the discrete least-squares (DLS) approximation of solutions of stochastic elliptic PDEs is improved by incorporating a reduced-basis method into the DLS framework. {In particular, we consider stochastic elliptic PDEs with an affine dependence on the random variable.}  The goal is to recover the entire solution map from {the} parameter space to the finite element space. To this end, first, a reduced-basis solution using {a} weak greedy algorithm is constructed, then a DLS approximation is determined by evaluating the reduced-basis approximation instead of the full finite element approximation. The main advantage of the new approach is that one only need apply the DLS operator to the coefficients of the reduced-basis expansion, resulting in huge savings in both the storage of the DLS coefficients and the online cost of evaluating the DLS approximation. In addition, the recently developed quasi-optimal polynomial space is also adopted in the new approach, resulting in superior convergence rates for a wider class of problems than previous analyzed. Numerical experiments are provided that illustrate the theoretical results.


\end{abstract}

\section{Introduction}

Mathematical models are used to understand and predict the behavior of complex systems arising in applications. Common input data for these type of models include forcing terms, boundary conditions, model coefficients, and the computational domain itself. Often, for any number of reasons there is a degree of uncertainty involved with these inputs. In order to obtain an accurate model one must incorporate such uncertainties into the governing equations and quantify their effect on model outputs of interest. In this paper, we focus on systems that can be modeled by elliptic partial differential equations (PDEs) with random input data {that has an affine dependence on the random variables}. In particular, we consider cases for which that data is parameterized, i.e., the random coefficients/fields in the PDEs are functions of a {\em finite} number of random parameters. This means that the PDE solution, denoted by $u(\bx, \bm y)$, can be viewed a function of an $N$-dimensional vector random parameters, denoted by $\bm y = (y_1, \ldots, y_N)^{\top}$. Our goal is to recover the entire solution map $\bm y \rightarrow u(\bx, \bm y)$ from the parameter space to the solution space. 

It is well known that commonly used Monte Carlo methods are not feasible options for this task because they can only be used to compute limited types of statistics. Sparse polynomial approximations \cite{Gunzburger:2014hi,nobile2008sparse,Cohen:2015jr,CDS2011,CCS2014,CCMFT2013}, including stochastic Galerkin, stochastic collocation, discrete least-squares (DLS), and compressive sensing methods, etc. These methods take advantages of the smoothness of the solution map $\bm y \rightarrow u(\bx, \bm y)$ to reduce the complexity of approximating that map in high-dimensional parameter space. The purpose of building sparse approximations (surrogates) is to enable fast evaluations of the surrogates, e.g., when conducting uncertainty quantification (UQ) tasks. However, existing sparse approximation techniques only focus on complexity reduction with respect to the parameter dependence, and largely ignore the huge cost (in evaluating the surrogates) arising from the finite element discretization. Specifically, denote by $J$ the degrees of freedom of the finite element discretization and by $M$ the dimension of the sparse polynomial space. To approximate the entire solution map $\bm y \rightarrow u(\bx, \bm y)$, we have to build a polynomial approximation for each of the $J$ finite element coefficients so that the final sparse approximation requires an $M \times J$ {\em dense} matrix to store all the coefficients. When $J$ is large, as it is in practical applications, the required storage may well not be affordable. Moreover, the complexity of each evaluation of the sparse approximation will be roughly of $\mathcal{O}(JM)$, which is considered as a perhaps prohibitive cost, given that computing accurate statistical information may require a very large number of such evaluations.

To overcome the above challenges, we propose to incorporate the well-studied reduced-basis technique \cite{BCDDPW2011,BBLMNP2010,HRS15,RHB2008} into the sparse polynomial approximation. Although we focus on improving discrete least-squares methods, our approach can also be potentially generalized to the other aforementioned approaches. The main idea is to construct a reduced-basis solution using a greedy algorithm \cite{BCDDPW2011} to reduce the number of coefficients that are dependent on the parameter vector $\bm y$ and then apply the DLS operator only to the coefficients of the reduced-basis approximation. For example, if the dimension $K$ of the reduced-basis space is such that $K \ll J$, then the DLS coefficients can be stored in an $M \times K$ matrix that is much smaller than the straightforward DLS case. Moreover, with respect to computational complexity, our approach requires only $\mathcal{O}(K(M+J))$ operations to perform matrix-vector productions for each evaluation of the DLS approximation, which is also significantly smaller to the $\mathcal{O}(MJ)$ operations required for the straightforward  DLS approximation.

The plan for the rest of the paper is as follows. In Section \ref{problemSetting}, we introduce the mathematical setting and assumptions needed throughout the rest of the paper. In Section \ref{Quasi}, we introduce quasi-optimal polynomial spaces and, in Section \ref{LeastSquares}, we discuss the formulation of the least-squares problem in Hilbert spaces using the quasi-optimal polynomial space introduced in Section \ref{Quasi}. Then, in Section \ref{ReducedBasis}, we introduce the reduced-basis method and discuss its incorporation into the least-squares framework. In Section \ref{numex}, we discuss the computational complexity of the discrete least-squares method as well as the reduced-basis method and provide the results of numerical experiments that illustrate our findings.

\section{Problem setting}\label{problemSetting}

Let $D$ denote a bounded Lipschitz domain in $\mathbb{R}^{d}$, $d \in \{1,2,3\}$, with boundary $\partial D$. We consider the parameterized elliptic partial differential problem
\begin{equation}\label{diffusion_deter}
\left\{\begin{aligned}
-\nabla \cdot \big(a({\bf x},\bm y)\nabla u({\bf x},\bm y)\big) &= f({\bf x}) \quad \forall ({\bf x},\bm y) \in D \times \Gamma\\
u({\bf x},\bm y) &= 0 \qquad \ \, \forall ({\bf x},\bm y)  \in \partial D \times \Gamma
\end{aligned}\right.
\end{equation}
for the unknown function $u({\bf x},\bm y)$, where $f({\bf x})$ and $a({\bf x},{\bm y})$ are given functions and ${\bm y}$ denotes a vector of parameters. In the stochastic setting, we have that ${\bm y}$ is a random variable distributed according to a joint probability density function (PDF) $\varphivarphi({\bm y})$. Although the case of a countably infinite number of random variables is of interest in some applications, here we assume, as is often the case, that the randomness present in a stochastic PDE can be approximated well in terms of a finite number of random variables. Therefore, we  assume that the parameterized diffusion coefficient $a({\bf x},\bm y)$ depends on a finite number $N$ of random variables denoted by the vector $\bm y = (y_{1},\dots, y_{N})^{\top} \in \Gamma\subset \mathbb{R}^N$, where $\Gamma$ denotes a parameter domain. 

Here, we specialize to the case for which the components of ${\bm y}$ are independent and identically distributed random variables so that $\Gamma$ is a hyper-rectangle in $\mathbb{R}^N$ and $\varphivarphi({\bm y}) = \Pi_{n =1}^{N}\varphivarphi_{n}(y_{n})$ is a product of $N$ one-dimensional PDFs. Without loss of generality, we can then assume that $\Gamma=[-1,1]^N$. We note that although the assumption that the random variables are i.i.d. may appear restrictive, in practice, a wide range of problems can still be addressed. As discussed in \cite{BNTT2014}, problems with non-independent random variables can be addressed via the introduction of auxiliary density functions. 

We make several assumptions. First, we assume that there exist constants $0 < a_{\min} < a_{\max} < \infty$ such that
\begin{equation}\label{assump12}
 a_{\min} < a({\bf x},{\bf y}) < a_{\max} \quad \mbox{$\forall\, {\bf x} \in D$,\,\, a.s. for \,${\bm y}\in\Gamma$}.
\end{equation}
Let $Y=L^{2}_{\varphivarphi}(\Gamma)$ denote the space of square integrable functions on $\Gamma$ with respect to the weight $\varphivarphi({\bm y})$. We also have the standard Sobolev space $X = H_{0}^{1}(D)$ equipped with the norm $\|v\|_{X}= (\int_D  \nabla v \cdot \nabla v d\bx)^{1/2}$; $X'=H^{-1}(D)$ denotes the corresponding dual space. Then, a weak formulation of \eqref{diffusion_deter} is, given $f\in X'$, to find $u\in X\otimes Y$ such that 
\begin{equation}
\label{diffusion_weak}
\begin{aligned}
\int_\Gamma  \int_D a({\bf x},{\bm y})\nabla& u({\bf x},{\bm y})\cdot \nabla v({\bf x},{\bm y})\varphivarphi(\bm y) d{\bf x}d{\bm y}
\\&
 =  \int_\Gamma \int_D  f({\bf x})  v({\bf x},{\bm y})\varphivarphi(\bm y) d{\bf x}d{\bm y}\ \ \ \forall\, v\in X\otimes Y,
\end{aligned}
\end{equation}
where $X\otimes Y := L_\varphivarphi^2(\Gamma;X) := \{ u : \int_\Gamma \|u(\cdot,{\bm y})\|_X^2\varphivarphi({\bm y})d{\bm y} < \infty\}$. By \eqref{assump12} and the Lax-Milgram theorem, there exists a unique solution of (\ref{diffusion_weak}) for any $f\in X'$ and that solution satisfies the bound
\begin{equation}
\|u\|_{X\otimes Y} \leq \frac1{a_{\min}} \|f\|_{X'}.
\end{equation}

Our final assumptions about the coefficient $a({\bf x},\bm y)$ is affine dependence on the random variables, i.e., it can be written in the form
\begin{equation}\label{KL}
a({\bf x},\bm y) = a_{0}({\bf x}) + \sum_{n=1}^{N}a_{k}({\bf x})y_{k}
\end{equation} 
for some $a_n({\bf x})$, $n\in\{0,\ldots,N\}$. A coefficient of this form could be a truncated Karhunen-Lo\`{e}ve (KL) expansion. The affine dependence is necessary to achieve satisfactory efficiency in constructing a reduced basis using greedy algorithms. {We also note that this assumption implies the complex continuation of $a({\bf x},\bm y)$, represented as the map $a({\bf x},\bm y):\mathbb{C}^N\to L^\infty ({D})$, is an $L^\infty({D})$-valued {\em holomorphic function} on $\mathbb{C}^N$. This allows for the use of the quasi-optimal polynomial space introduced in Section \ref{Quasi}.} 

For spatial discretization, we use standard finite element methods. Let $X^{\rm fe}_{h}\subset X$ denote a standard finite element space of dimension $J$ and let  $\{\phi_j({\bf x})\}_{j=1}^J$ denote a basis for $X^{\rm fe}_{h}$ that consists of piecewise-continuous  polynomials defined with respect to a regular triangulation $\mathcal{T}_{h}$ of $D$, where $h > 0$ denotes the maximum mesh spacing. For any $\bm y \in \Gamma$, the finite element approximation $u_h({\bf x}, \bm y) \in X_h^{\rm fe}$ is determined by solving 
\begin{equation}
\label{BillinearWeak}
A(u_h(\bm y),v;\bm y) = (f,v) \ \ \ \forall v \in X_h^{\rm fe},
\end{equation}
where $A(u,v) := \int_D a({\bf x}, \bm y) \nabla u({\bf x},\bm y) \cdot \nabla v({\bf x},\bm y) d{\bf x}$ is the bilinear form corresponding to the operator in \eqref{diffusion_weak} and $(\cdot,\cdot)$ denotes the $L^2(D)$ inner product. Then, for any $\bm y \in \Gamma $, we have, for sufficient small $h$ and for a constant $C_S$ whose value is independent of $h$ and $\bm y$, the error estimate 
\begin{equation}\label{fees}
\|u - u_{h}\|_{X} \leq C_S h^{\alpha}.
\end{equation}
The convergence rate $\alpha$ depends on the spatial regularity of $u$ and the degree of the polynomial used. For a detailed treatment of finite element error analyses, see, e.g., \cite{BS2008}.

\section{Discrete least-squares approximation}\label{DLS}

In this section, we recall the formulation and theoretical results about random discrete $L^{2}$ approximations of solutions of the parameterized PDE problem \eqref{diffusion_deter}. This presentation is brief and only discusses the DLS approximation for Hilbert-valued functions. For a more general and comprehensive analysis, see, e.g., \cite{CCMFT2013,MFST2013}. For the sake of further simplifying the exposition, we assume that the random variables are uniformly distributed on {$[-1,1]$} so that $\varphivarphi_n=\frac12$ for $n = 1, \ldots, N$.

\subsection{Quasi-optimal polynomial spaces}\label{Quasi}

The first step towards building a DLS approximation is to choose an appropriate polynomial space in $L_\rho(\Gamma)$. Because we assume a uniform measure, we use Legendre polynomials that are orthogonal with respect to this measure. Letting $\bm \nu = (\nu_{1},\nu_{2}, \dots,\nu_{N}) \in \mathbb{N}_{0}^{N}$ denote a multi-index, the multidimensional Legendre polynomials are denoted by $L_{\bm \nu}(\bm y) = \prod_{i = 1}^{N} L_{\nu_{n}}(y_{n})$, where $L_{\nu_{n}}(y_{n})$ denotes the $L^{2}$-normalized one-dimensional Legendre polynomials \cite{2013JSV...332.4403B}.

In the construction of polynomial approximations with respect to parameter dependences, one wishes to select a multi-index set $\Lambda_M \subset \{ \bm \nu = (\nu_{1},\nu_{2}, \dots,\nu_{N})\,:\, \nu_n \in \mathbb{N}_0\}$ such that the corresponding polynomial space $\mbox{\em span}\{L_{\bm \nu}(\bm y), \bm \nu \in \Lambda_M\}$ yields maximal accuracy for a given dimension $\mmm$. To achieve this, there are two approaches that have been extensively studied \cite{CCS2014,CDS2010,CDS2011}. The first approach is known as {\em best $M$-term approximation}. Using a truncated Legendre expansion and using the triangle inequality, we can express the error of the approximation in the form 
\begin{equation}
\label{bestErr}
 \Big\| u(\bm y) - \sum\limits_{{\bm \nu} \in \Lambda_M}  c_{\bm \nu}L_{\bm \nu}(\bm y) \Big\|_{X\otimes Y}\!\! \leq \sum\limits_{{\bm \nu} \in \Lambda_M} \|c_{\bm \nu}\|_{X},
\end{equation} 
where $\Lambda_{M}$ is chosen such that the error (\ref{bestErr}) is minimized. This means that the indices $\bm \nu \in \Lambda_M$ correspond to the $\mmm$ largest values of $\|c_{\bm \nu}\|_{X}$. However, in practice, finding the best index set and polynomial space is an infeasible task because it requires computation of {\em all} the coefficients $c_{\bm \nu}$. 
 
An alternative approach that tends to be less computationally intensive is referred to as {\em quasi-optimal polynomial approximation} \cite{BNTT2014,Tran2017}. Rather than explicitly computing the coefficients in order to evaluate $\|c_{\bm \nu}\|_{X}$, we instead compute sharp estimates for $\|c_{\bm \nu}\|_{X}$ and use these to determine a quasi-optimal index set $\Lambda_{M}$. It has been shown that this method can achieve convergence rates similar to those of the best $\mmm$-term approximation.

In order to establish a bound on the coefficients of the Legendre expansion, we need the following definition concerning uniform ellipticity in polyellipses.

\begin{definition}
\label{cext}
For $0 < \delta < a_{\min}$ and $\bm \rhorho$ denoting the sequence $\{\rhorho_{i}\}_{i=1}^N$ with $\rhorho_{i} >1$  $\forall i$, we say the random field $a(\bf x,\cdot)$ satisfies the {\em$(\delta,\bm \rhorho)$-polyellipse uniform ellipticity} assumption if it holds that 
$$
\mathcal{R}(a(\bf x,\bm z)) \geq \delta
$$
for all $\bx \in D$ and $\bm z = \{z_{i}\}_{i=1}^N$ contained in the polyellipse 
$$
\mathcal{E} = \bigotimes_{1 \leq i \leq N} \left\{z_{i} \in \mathbb{C}: \mathcal{R}(z_{i}) = \frac{\rhorho_{i} + \rhorho_{i}^{-1}}{2} \cos(\theta), \mathcal{I}(z_{i}) = \frac{\rhorho_{i} - \rhorho_{i}^{-1}}{2} \sin(\theta), \theta \in [0,2\pi) \right\}.
$$
\end{definition}
It has been shown \cite{Tran2017}, for any diffusion coefficient $a({\bf x},\bm y)$ satisfying the coercivity assumption in \eqref{assump12} and having the holomorphic parameter dependence, there always exists one $\bm \rhorho$ for which this property is satisfied. Now, using this regularity condition, the holomorphy of the solution with respect to the random parameters follows and the bound on the coefficients of the $L^{2}$-normalized Legendre expansion 
\begin{equation}\label{legendre_est2}
\|c_{\bm \nu}\|_X \le C_{{\bm \rhorho},\delta} {\bm \rhorho}^{-{\bm \nu}}\prod_{i=1}^N\sqrt{2\nu_i+1}
\end{equation}
 holds, where $C_{{\bm \rhorho},\delta}  = \frac{||f||_{V'}}{\delta}\prod_{i = 1}^{N} \frac{\ell(\mathcal{E}_{\rhorho_{i}})}{4(\rhorho_{i}-1)}$ with $\ell(\mathcal{E}_{\rhorho_{i}})$ denoting the perimeter of the ellipse $\mathcal{E}_{\rhorho_{i}}$. Note that Definition \ref{cext} holds for an infinite combination of $(\delta,\bm \rhorho)$ that we denoted by $\bm A \bm d$. For a given $\bm \nu$, the best coefficient bound will then be given by
$$
\|c_{\bm \nu}\|_X \le \inf_{(\delta,\bm \rhorho)\in \bm A \bm d}C_{{\bm \rhorho},\delta} {\bm \rhorho}^{-{\bm \nu}}\prod_{i=1}^N\sqrt{2\nu_i+1}.
$$
Solving this minimization problem is in general computationally infeasible. However, in case the basis functions $a_k$ have non-overlapping supports, $\bm \rhorho$ can be determined easily \cite{BNTT2014}. Problems with both overlapping support and nonoverlapping support are explored further in Section \ref{numex}. We can now state an asymptotic bound for the quasi-optimal $\mmm$-term approximation as follows:
{
\begin{proposition}
	\label{leg_bound}
	Consider the Legendre series $\sum_{{\bm \nu}\in \Lambda} c_{\bm \nu}L_{\bm \nu}$ for $u$. Assume that \eqref{legendre_est2}
	holds for all ${\bm \nu}\in\Lambda$. Let $\log(\rhorho_{n}) = \lambda_{n}$
	and $\Lambda_M$ denote the set of indices corresponding to the $\mmm$ largest bounds in \eqref{legendre_est2} determined by  
		\begin{equation}
		\label{logCorrection}
		\Lambda_M = {\Bigg\{{\bm \nu} \in \Lambda :  \sum\limits_{n=1}^N (2\lambda_n \nu_{n} - \log (2\nu_{n} + 1)) \leq j\Bigg\}}
		\end{equation} 
	for a given $j\in \mathbb{N}$. Then, for any $0<\mu<1$, it follows that
	\begin{align}
	\label{theorem:est2}
	\Big\| u - \sum\limits_{{\bm \nu} \in \Lambda_M}  c_{\bm \nu}L_{\bm \nu} \Big\|_{X\otimes Y}^2\!\! \le C_{{\bm \rhorho},\delta}^2 C_u(\mu) \mmm\exp\! \bigg(\! - 2 \Big({\frac{\mmm N! \prod_{n=1}^{N}\lambda_n}{(1+\mu)}}\Big)^{\frac{1}{N}} \bigg)
	\end{align}
	as $j$ (and therefore M) {goes to $\infty$}.  Here, 
	$
	C_u(\mu) = (4e + 4\mu e- 2) \frac{e}{e-1}.
	$
\end{proposition}
}
{ The index set should be chosen such that its size $M$ allows for a specified level of accuracy to be reached using estimate \eqref{theorem:est2}.  } The value $0 < \mu < 1$ is related to the cardinality of our polynomial approximation, and decreases towards $0$ as the cardinality of our polynomial approximation increases. A sharp mathematical formula for $\mu$, given $\mmm$, is currently an open problem, though it has been shown that even for a moderate value of $\mu$ one can still obtain a strong rate of convergence.  For full details see \cite[Section 4]{Tran2017}. 

\subsection{Discrete least-squares approximation in quasi-optimal spaces}
\label{LeastSquares}

Here we introduce the DLS method for approximating solutions of parametric PDE in \eqref{diffusion_deter} in the quasi-optimal polynomial space discussed above. This presentation is brief and only discusses the least-squares approximation for Hilbert-valued functions. For a more general and deeper analysis, see, e.g., \cite{CCMFT2013,MFST2013}.

Let $Y_{\Lambda_M}$ denote an $M$-dimensional quasi-optimal subspace in $Y$. We intend to build a DLS approximation in $X\otimes Y_{\Lambda_M}$ of the solution $u({\bf x},{\bm y})\in X\otimes Y$ by the orthogonal projection, i.e.,
$$
P_{M}[u] := \argmin_{v\in X \otimes Y_{\Lambda_M}} \|u - v\|_{X\otimes Y}.
$$
Letting $\bm{\ell}_{M}(\bm y) := (\ell_{1}({\bm y}),\ldots, \ell_{M}(\bm y))^{\top}$ denote the vector of { re-indexed} Legendre basis functions $\{L_{\bm \nu}(\bm y): \bm \nu \in \Lambda_M\}$ for the subspace $Y_{\Lambda_M}$, we have that
$$
  P_{M}[u] = \sum_{m=1}^M c_m({\bf x}) \ell_m({\bm y})
  \quad\mbox{with}\quad 
  c_m({\bf x}) = \big\langle  u({\bf x},\cdot), \ell_m(\cdot) \big\rangle_Y \quad\mbox{for $m=1,\ldots,M$},
$$
where $\langle\cdot,\cdot\rangle_Y$ denotes the inner product on $Y$.

In general, we do not have available the solution of the PDE for all ${\bm y}\in\Gamma$, but only at a set of points $\{{\bm y}_i\}_{i=1}^S$, where ${\bm y}_i\in\Gamma$ are i.i.d. random variables distributed according to some distribution. We then consider the discrete (with respect to the ${\bm y}$ dependence) least-squares problem
\begin{equation}
\label{discreteLQ}
P_{M,S}[u]:= \argmin_{v\in X \otimes Y_{\Lambda_M}} \sum_{i = 1}^{S} \|u({\bf x},{\bm y}_{i}) - v({\bf x},{\bm y}_{i})\|^2_X
\end{equation}
that has a unique solution as long as $M \leq S$.


In practice, we do not have access to the exact solution $u({\bf x},{\bm y}_i)$ for ${\bm y}_i\in\Gamma$, so that we apply the DLS operator to the finite element solution
$u_h \in X_h^{\rm fe} \otimes Y$ and obtain the $L^2$ projection $P_{M,S}[u_h]$ in the subspace $X_h^{\rm fe} \otimes Y_{\Lambda_M}$ which has the form

\begin{equation}\label{DLS_uk}
P_{M,S}[u_h] = \sum_{m = 1}^M \sum_{j = 1}^J c_{mj}  \phi_j({\bf x})\ell_m(\bm y),
\end{equation}
{recalling that $\phi_j({\bf x})$ are the finite element basis functions introduced in section \ref{problemSetting}.} {Letting $[{\bm\Phi}]_{ij} = \ell_j(\bm y_i)$}, the coefficients $\{c_{mj}\}_{m=1,j=1}^{M,J}$ are the solution of the following linear system:
\begin{equation}\label{lspde}
({\bm\Phi}^{\top}{\bm\Phi}) \mathbf{C} = {\bm\Phi}^{\top} \mathbf{U},
\end{equation}
where $[\mathbf{C}]_{mj} = c_{mj}$ for $m = 1, \ldots, M$, $j = 1, \ldots, J$ and $[\mathbf{U}]_{ij} = u_h({\bf x}_j, \bm y_i)$ for $i = 1, \ldots, S$, $j =1, \ldots, J$.

\section{Improved DLS methods based on reduced-basis solutions}
\label{ReducedBasis}

The main purpose of building DLS approximations is to reduce the cost of obtaining approximate solutions of the PDE problem \eqref{diffusion_deter} at a large set of samples in $\Gamma$, i.e., reducing the online cost. We observe that the need to reduce costs is only necessary when the finite element degrees of freedom $J$ is extremely large. Otherwise, for a small $J$, a classic finite element solver will be efficient enough to be used as an online solver. However, for a very large $J$, we can see from \eqref{lspde} that the coefficient $\mathbf C$, which is an $M \times J$ {\em dense} matrix, may require an unaffordable amount of storage. Moreover, the complexity of each evaluation of the DLS approximation would be of $\mathcal{O}(JM)$. To avoid such inefficiencies in both storage and computation, we propose to develop a new DLS method based on reduced-basis approximations of solutions of \eqref{diffusion_deter}. A brief overview of a reduced-basis method is given in Section \ref{sec:RB} and our approach is introduced in Section \ref{sec:RB_DLS}.

\subsection{Reduced-basis methods}\label{sec:RB}

We briefly recall the reduced-basis technique; for a more in depth discussions about reduced-basis methods, see \cite{RHB2008,BBLMNP2010}. The main idea of reduced-basis methods is to collect a set of deterministic solutions of the stochastic problem in \eqref{diffusion_deter} at a subset of the most representative samples in $\Gamma$, then uses these solutions as a basis to approximate solutions at other points in $\Gamma$ through Galerkin projection. Specifically, when we have a subset, denoted by $\Xi_\kkk := \{\bm y_i\}_{i=0}^{\kkk-1}$, consisting of $\kkk$ representative  samples, we can then solve the finite element system in \eqref{BillinearWeak} $\kkk$ times to obtain the set of $\kkk$ solutions 
$\{u^i_h({\bf x})=u_h({\bf x},\bm y_i)\}_{i=0}^{\kkk-1}$. Using these solutions (snapshots), we can define a $\kkk$-dimensional reduced space 
\[
X^{\rm rb}_\kkk={\rm span}\{u^i_h(\bx)\}_{i=0}^{\kkk-1}\subset X^{\rm fe}_h.
\]
Then, we can construct a reduced-basis approximation $u_{h,\kkk}$ by projecting $u_h$ into $X^{\rm rb}_\kkk$, i.e., seeking
\begin{equation}\label{rbs0}
u_{h,\kkk}(\bx,\vy) = \sum_{i=0}^{\kkk-1} w^k_i(\vy)\, \xi_i(\bx)\in X^{\rm rb}_\kkk,
\end{equation}
satisfying
\begin{equation}\label{nlwprb}
A(u_{h,\kkk}(\bm y),v;{\bm y}) =(f,v)
\quad\forall\,v\in X^{\rm rb}_\kkk,
\end{equation}
where $A(u,v;\vy)$ is the bilinear form defined in \eqref{BillinearWeak} and $\{\xi_i\}_{i=0}^{\kkk-1}$ is the orthogonalized reduced basis of $X^{\rm rb}_\kkk$. Note that, for each $\vy\in\Gamma$, the equation in \eqref{nlwprb} is equivalent to a linear system of $\kkk$ algebraic equations for the coefficients $\{w_i^k(\vy)\}_{i=0}^{\kkk-1}$ in \eqref{rbs0}. In this way, when $\kkk \ll J$, the computational cost of  approximating $u(\bx,\vy)$ for each $\bm y \in \Gamma$ is significantly reduced from solving a $J\times J$ linear system to solving a $\kkk \times \kkk$ linear system. 

Now the question is how does one determine a good set of $\kkk$ samples $\Xi_\kkk := \{\bm y_i\}_{i=0}^{\kkk-1}$? Suppose one has in hand a set of samples $\Xi_k$; one could start with $k=0$, i.e., a single sample chosen at random or at the center of $\Gamma$. Then given the current set $\Xi_k$ of samples, how does one find the next sample $\bm y_k$ in an effective and efficient way so as to improve the accuracy of the reduced-basis solution. The ideal choice is to use the {\em greedy algorithm} \cite{BCDDPW2011}, i.e., find the next sample $\bm y_k$ by solving the optimization problem
\begin{equation}\label{greedygreedy}
\bm y_{k} =  \argsup_{\bm y\in\Gamma}\|{ u_{h,k}(\cdot,\vy) - u_h(\cdot,\vy) }\|_{X}, 
\end{equation}
i.e., locating the point $\bm y_{k}\in\Gamma$ at which the error between the current reduced-basis approximation and the finite element approximation is the largest. Unfortunately, solving the optimization problem \eqref{greedygreedy} is not practical because it requires full information about the exact finite element solution $u_h(\bx, \vy)$. To circumvent this issue, a variant of the greedy strategy \eqref{greedygreedy}, i.e., the {\em weak greedy algorithm}, has been shown to be computationally feasible in the context of solving parameterized PDEs \cite{RHB2008}. The key idea of the weak greedy strategy is to find an accurate and computationally efficient surrogate of the error $u_{h,k}(\cdot,\vy) - u_h(\cdot,\vy) $, and replace the true error in \eqref{greedygreedy} with the surrogate to solve the optimization problem. To this end, we use the Galerkin residual as the surrogate to implement the weak greedy algorithm. Letting $e_{h,k}(\bx,\vy):=u_h(\bx,\vy) - u_{h,k}(\bx,\vy)\in X^{\rm fe}_h$, we then have that, for any $\vy\in\Gamma$,
\begin{equation}\label{resid3}
 R(v;\bm y ) := A\big(e_{h,k}(\cdot,\vy),v; \vy\big) = (f,v) - A\big(u_{h,k}(\cdot,\vy),v; \vy\big)\quad \forall\, v\in X^{\rm fe}_h.
\end{equation}
Thanks to Riesz representation, we have $(\hat{e}_{h,k},{e}_{h,k})_X = R({e}_{h,k};\bm y)$, such that
\[
\|e_{h,k}(\bx,\vy)\|_X=\|u_h(\bx,\vy) - u_{h,k}(\bx,\vy)\|_X \le \frac{1}{\alpha_{\rm LB}(\bm y)}\|\hat{e}_{h,k}\|_X,
\]
where $\alpha_{\rm LB}(\bm y) = \min_{\bx \in D} a(\bx, \vy)$.  Thus, we can replace $e_{h,k}$ with $\hat{e}_{h,k}$. In this effort, we also have to replace the search over all $\vy\in\Gamma$ by a search over a discrete training set; for solutions manifolds that are sufficiently smooth, this step does not introduce unmanageable errors. Specifically, the construction of the reduced-basis method begins by choosing a training set $\Xi_{\rm train}$ of $S_{\rm train}$ points in $\Gamma$; these points could be chosen randomly according to the joint PDF $\rho(\vy)$ associated with the random parameters $\vy\in\Gamma$ or could be chosen deterministically. Then, the optimization problem in \eqref{greedygreedy} is solved within the training set, i.e., $\bm y_k$ is generated by 
\begin{equation}\label{greedy3}
    \vy_{k} = \argsup_{\vy\in\Xi_{\rm train}}\Big\{\frac1{\alpha_{\rm LB}(\vy)} \| {\hat e}_{h,k}(\cdot,\vy) \|_X\Big\}.
\end{equation}
{The term $\alpha_{\rm LB}(\bm y)$ will be calculated over the training set using the Successive Constraint Method outlined in \cite{HRS15}}

Due to the affine property of the coefficient in \eqref{KL}, the residual $R(v;\bm y)$ in \eqref{resid3} can be decomposed as
\[
\begin{aligned}
R(v;\bm y) & = (f,v) - A_0\big(u_{h,k}(\cdot,\vy),v\big)  + \sum_{n=1}^{N}A_n\big(u_{h,k}(\cdot,\vy),v;\bm y\big) y_{n}\\
& = (f,v) - \sum_{i=0}^{k-1} \Big(A_0\big(\xi_i,v\big)  + \sum_{n=1}^{N}A_n\big(\xi_i,v\big) y_{n}\Big) w_i^k(\bm y),
\end{aligned}
\]
for all $v \in X_h^{\rm fe}$. Due to the linearity of the above representation, we can determine $\hat e_{h,k}(\bx,\vy)\in X^{\rm fe}_h$ efficiently by solving the following set of problems offline:
\begin{equation}\label{resid7}
\left\{\begin{aligned}
   &(\hat e^f,v)_X = (f,v)\quad \quad\;\,\forall\,v\in X^{\rm fe}_h,\\
   &(\hat e_{n,i},v)_X =  A_n(\xi_i,v)\quad\forall\,v\in X^{\rm fe}_h \quad \mbox{for $i=0,\ldots,k-1$ and $n=0,\ldots,N$},
\end{aligned}\right.
\end{equation}
such that $\hat{e}_{h,k}$ can be computed very efficiently by
$$
{\hat e}_{h,k}(\bx,\vy) = \hat e^f(\bx) - 
\sum_{i=0}^{k-1} w_i^k(\vy) \Big[\hat e_{0,i}(\bx)
+\sum_{n=1}^N y_n \hat e_{n,i}(\bx)
\Big]\in X^{\rm fe}_h.
$$
{We note that the quantities $\hat e^f(\bx)$, $\hat e_{0,i}(\bx)$, and $\hat e_{n,i}(\bx)$ are independent of $\bm y$ and can therefore be stored in an offline phase.}

To terminate the greedy procedure, we can preset some error tolerance $\varepsilon_{\rm tol}$ and end the algorithm when the approximation error is judged to be sufficiently small, i.e., $\{\|\hat{e}_{h,k}/\alpha_{\rm LB}(\bm y)\|_X\leq \varepsilon_{\rm rb}, \forall \bm y \in \Xi_{\text{train}}\}$. In addition, how the points should be selected and the size of the training set is extremely problem dependent. In practice there are two approaches commonly used to construct $\Xi_{\text{train}}$. The first is an adaptive approach which starts with a small number of sample points and then greedily enriches the sample space based on these initial points; see \cite{HSZ2011}. The other method is to randomly sample the parameter space $\Gamma$ according to the probability distribution associated with the problem.

\subsection{Reduced-basis discrete least-squares (RB-DLS) approximation} \label{sec:RB_DLS}

As already mentioned, the goal of using reduced-basis approximations in the least-squares setting is to reduce the online cost, i.e., the cost of evaluating the final DLS approximation. Letting $K$ denote the final value of $k$ upon termination of the greedy algorithm, we observe that the parameter dependence of the reduced-basis solution in \eqref{rbs0} only appears in the coefficients $\bm w_\kkk(\bm y) := (w^\kkk_1(\bm y), \ldots, w_\kkk^\kkk(\bm y))^{\top}$ which is vector of size $K \ll J$. Thus, instead of applying the DLS operator $P_{M,S}[\cdot]$ to $u_h$, we apply it to the reduced-basis solution $u_{h,\kkk}$, i.e., 
\begin{equation}\label{RB_DLS}
P_{M,S}[u_{h,\kkk}] = \sum_{k = 0}^\kkk \left(\sum_{m = 1}^M  c_{m k}^{\rm rb} \, \ell_m(\bm y)\right) \xi_k({\bf x})
\end{equation}
which is equivalent to approximating the coefficient vector $\bm w_K(\bm y)$ using the DLS method. The algebraic formulation for solving the coefficients $c_{mk}^{\rm rb}$ is
\begin{equation}\label{lspde1}
({\bm\Phi}^{\top}{\bm\Phi}) \mathbf{C}^{\rm rb} = {\bm\Phi}^{\top} \mathbf{W},
\end{equation}
where $[\mathbf{C^{\rm rb}]}_{mk} = c_{mk}^{\rm rb}$ for $m = 1, \ldots, M$, $k = 0, \ldots, \kkk$ and $[\mathbf{W}]_{ik} = w_i^K(\bm y_k)$ for $i = 1, \ldots, S$, $k =0, \ldots, K$. Then, for each new sample $\bm y \in \Gamma$, the evaluation of the RB-DLS approximation in $P_{M,S}[u_{h,\kkk}]$ can be conducted by
\begin{equation}\label{uRBDLS}
 \bm u^{\rm RB-DLS} =  \mathbf{V}(\mathbf{C^{\rm rb}})^{\top} \bm{\ell}_{M}(\bm y),
\end{equation}
where $\bm u^{\rm RB-DLS}:= (P_{M,S}[u_{h,\kkk}](\bx_1), \ldots, P_{M,S}[u_{h,\kkk}](\bx_J))^{\top}$ and $\mathbf{V} = (\xi_0, \ldots, \xi_\kkk)$ is the reduced-basis matrix. In comparison, evaluating the classic DLS approximation $P_{M,S}[u_h]$ in \eqref{DLS_uk} has to be done by
\begin{equation}\label{uDLS}
\bm u^{\rm DLS} = \mathbf{C}^{\top} \bm{\ell}_{M}(\bm y),
\end{equation}
where $\bm u^{\rm DLS}:= (P_{M,S}[u_{h}](\bx_1), \ldots, P_{M,S}[u_{h}](\bx_J))^{\top}$ and $\mathbf{C}$ is given in \eqref{lspde}. The advantages of \eqref{uRBDLS} over \eqref{uDLS} can be seen from two aspects. In terms of storage, \eqref{uRBDLS} only requires storage for a $J\times \kkk$ matrix $\mathbf{V}$ and an $M \times \kkk$ matrix $\mathbf{C}^{\rm rb}$, whereas \eqref{uDLS} requires storage for an $M \times J$ matrix $\mathbf{C}$. Thus, when $\kkk$ is small, \eqref{uRBDLS} requires much less storage than \eqref{uRBDLS}. In other words, the matrix $\mathbf{C^{\rm rb}}\mathbf{V}^{\top}$ can be viewed as a {\em low-rank} (i.e., rank $\kkk$) approximation of the matrix $\mathbf{C}$. In terms of computation, for each $\bm y \in \Gamma$, \eqref{uRBDLS} requires $\mathcal{O}(\kkk(M+J))$ operations to perform matrix-vector products, whereas \eqref{uDLS} requires $\mathcal{O}(MJ)$ operations. This is another significant savings achieved by using our approach.

The total error of the RB-DLS approximation $P_{M,S}[u_{h,\kkk}]$ can be split into the sums of the finite element discretization error, the reduced-basis error, and the DLS projection error, i.e.,
\begin{equation}
\label{conv_rates}
\begin{aligned}
& \mathbb{E}\left[\|u - P_{M,S}[u_{h,\kkk}]\|_{X}^{2}\right]\\  \leq &  {\underbrace{ \mathbb{E}\left[\|u - u_{h}\|_{X}^{2}\right]}_{e_{\rm I}}  +  \underbrace{\mathbb{E}\left[\|u_{h} - u_{h,\kkk}\|_{X}^{2}\right]}_{e_{\rm II}}+ \underbrace{\mathbb{E}\left[\|u_{h,\kkk} - P_{M,S}[u_{h,\kkk}]\|_{X}^{2}\right]}_{e_{\rm III}}.}
\end{aligned}
\end{equation}
The first error $e_{\rm I}$ is easy to control/balance based on the classic finite element error analysis. The second error $e_{\rm II}$ is essentially controlled by the Kolmogorov width associated with the solution manifold in the finite element space $X_h^{\rm fe}$. In the recent works \cite{BC2015, CD2015},  it has been shown for the class of problems dealt with in this paper that the Kolmogorov width will decay {at least algebraically} ; this gives us hope that our reduced-basis method will be successful. {In practice it is difficult to determine the decay of the RB error a-priori, hence we use the {\em a-posteriori} estimate to balance the second error $e_{\rm II}$ by adjusting the threshold $\varepsilon_{\rm tol}$}.
 
Defining the term
\[
\KKK(\Lambda_{M}):=  \sum_{\bm \nu \in \Lambda_{{M}}}||L_{\bm \nu}||^{2}_{L^{\infty}},
\]
the third error $e_{\rm III}$ can be bounded {for any $r > 0$} by \cite{CDL2013}
\[
e_{\rm III} \le (1 + \beta(S))e_{M}(u)^{2} + 8 H^{2}S^{-r}
\]
as long as the number of sample points satisfies
\[
\frac{S}{\ln(S)} \geq \frac{\KKK(\Lambda_{M})}{\kappa},\ \ \ \kappa := \frac{1 - \ln(2)}{2 + 2r},
\]
where  $\beta(S) \rightarrow 0$ as $S \rightarrow +\infty$, $H$ is the uniform upper bound of $u$, and $e_M(u)$ is the error in the best $M$-term approximation of $u_{h,K}$. {The term $r$ is related to the stability of the least squares system, full details can be found in  \cite{CDL2013}}. { As shown in \cite{CCMFT2013}, when using Legendre polynomials, the quantity $Z(\Lambda_M)$ can be bounded by $M \le Z(\Lambda_M) \le M^2$, when $\Lambda_M$ is a lower set. When using Chebyshev polynomials, a better bound can be obtained for lower sets, i.e., $\KKK(\Lambda_{M}) \leq  \min( M^{\log 3/\log 2},2^{N} M)$. Even though the polynomial space used in this work is not a lower set, the quasi-optimal polynomial space can be covered with a lower set which is only slightly larger allowing us to effectively use these bounds. As will be illustrated in section \ref{numex} to maintain stability of the DLS method using the quasi-optimal polynomial space sample points only on the order of $3M$ are required. Therefore the slight theoretical oversampling due to the lack of the quasi-optimal polynomial space not being a lower set does not have an impact in practice. }
 With the use of the quasi-optimal error estimate in Proposition \ref{leg_bound}, the error $e_{\rm III}$ can then be bounded by 
\begin{equation}\label{err3}
e_{\rm III} \leq (S^{\frac{N}{N+1}}C_{N} + 8H^{2})\exp(-(C_{e}S)^{\frac{1}{N+1}}),
\end{equation}
where the constants $C_e$ and $C_N$ are given by
\begin{equation}
\begin{aligned}
C_{e} &= \left(\frac{\tau 2^{N-1} N! \prod_{i=1}^{N} \lambda_{i} }{3^{N}(1+\mu)}\right) \\ 
C_{N} &=  ((1 + \beta(S))C_{{\bm \rhorho},\delta}C_{u}(\mu)\frac{\tau^{\frac{N}{N+1}}}{3 * 2^{2N}}\left(\frac{(1+\mu)}{N! \prod_{i=1}^{N} \lambda_{i}}\right)^{\frac{1}{N+1}}
\end{aligned}
\end{equation}
with $\tau = \frac{1 - \ln 2}{2}$ and the constants $\mu, C_{\bm \varphi, \delta}, C_{\mu}$ are defined in Section \ref{Quasi}.

The error $e_{\rm III}$ can be balanced by constructing an appropriate quasi-optimal subspace. To do so, we will be required to determine the weights $\bm\lambda = (\lambda_1, \ldots, \lambda_N)$ in \eqref{theorem:est2}. We note that it is the ratio of the weights which will determine our polynomial space, the magnitude will simply dictate the pace at which the error decays. It it only possible to analytically construct the weights in the case where the functions $\{a_n(\bx)\}_{n=0}^N$ in \eqref{KL} do not have overlapping supports, e.g., the inclusion problem investigated in \cite{BNTT2014}; otherwise the weights must be determined numerically. On the other hand, the optimal weights require the solution of a nonlinear optimization problem in the $N$-dimensional parameter space, which is also not feasible in practice. Hence,  we instead follow a procedure of one-dimensional analyses {as} done in \cite{BNTT2012,FTW2008}. We consider the subset $\mathbb{U} = \{ \bm \nu \in \mathbb{N}_{0}^{N} :  \nu_{i} = 0 \,\, \text{if} \,\, i \neq n, \nu_{n} = 0,1,2,. . . \}$. Then, according to the decay rates established in the previous section, $|c_{\bm \nu}| \sim e^{- \bm\lambda_{n}}$, so the rate $\bm\lambda_{n}$ can then be estimated through linear {regression} of the quantities $\ln|c_{\bm \nu}|$. {Now recalling definition 1 for any $0 \leq \delta \leq 1$ it must hold that $\mathcal{R}(a(\bf x,\bm z)) \geq \delta$ for all $\bx \in D$ and $\bm z = \{z_{i}\}_{i=1}^N$ in the polyellipse $\mathcal{E}$ which is determined by the weights $\bm\lambda_{n}$. In order to ensure this holds we scale our weights by an appropriate constant. Even though this may not result in an optimal estimate, it will still manage to capture any anisotropic behavior present in the problem.}

\section{Numerical experiments}
\label{numex}
In this section, we illustrate the convergence as well as the computational efficiency of the DLS-RB method. All calculations in this section are effected using the FEniCS \cite{Alnaes2012a} (\url{http://fenicsproject.org/}) and Rbnics \cite{HRS15} \url{http://mathlab.sissa.it/rbnics} software suites. 
 We will use the same problem formulation utilized in \cite{Chen2014}. Consider the stochastic elliptic problem \eqref{diffusion_weak} with $D = [0,1]^{2}$, the forcing term $f =  1$, and the finite element discretization with fixed $h = \frac{1}{256}$. We take the coefficient $a(\mathbf x,\bm y)$ to be a random field with expectation and correlation given as %
 \begin{equation}
 \mathbb{E}[a]({\bf x}) = c\; \ \text{for a fixed $c > 0$}\;\; \text{and} \;\; \mathbb{C}ov[a](\mathbf x,\mathbf x') = \exp \left(-\frac{(\mathbf x - \mathbf x')^{2}}{L^{2}}\right),
 \end{equation}
 where $L$ is the correlation length.  This field can be represented by the following Fourier-type expansion 
 \begin{equation}\label{klexp}
 a(\mathbf x,\textbf{y}) = \frac{1}{100}\Bigg\{c +\bigg(\frac{\sqrt{\pi}L}{2}\bigg)^{\frac{1}{2}}y_{1} + \sum_{n=1}^{\infty}\sqrt{\xi_{n}}\Big(\sin(n\pi x_{1})y_{2n} + \cos(n\pi x_{1})y_{2n + 1}\Big)\Bigg\},
 \end{equation}
 where the uncorrelated random variables $y_{n}$ have zero mean and unit variance, and the eigenvalues are equal to
 \begin{equation}
 \sqrt{\xi_{n}} = (\sqrt{\pi}L)^{\frac{1}{2}}\exp\Bigg(-\frac{(n\pi L)^{2}}{8}\Bigg) \;\;\text{for}\;\; n \geq 1.
 \end{equation}
Here, we take $c = 4$, $L = \frac{1}{8}$ and only retain the first 5 random variables $\bm y = (y_{1}, . . .,y_{5})$ in the expansion \eqref{klexp}. Even though the independence of the five random variables is only valid in the case of Gaussian distribution, we assume $(y_{1}, . . .,y_{5})$ are independent uniformly distributed random variables in $ \Gamma = [-1,1]^{5}$. The weights for the quasi-optimal subspace are found to be $\bm\lambda \approx (
0.68, 0.66, 0.98,1.37, 0.49)$ after rescaling. In order to measure the error in our examples we will consider the quantity of interest
 \begin{equation}
 Q(u) = \frac{1}{|D|}\int_{D}ud\mathbf x,
 \end{equation}
 and examine the behavior of our algorithm in the norm
 \begin{equation}
 \mathbb{E}\big(\|Q(u_{h}) - Q(P_{M,S}[u_{h,k}](\bm y)) \|_{\infty}\big) \approx \max_{\bm y \in \Xi_{\rm test}} \left|Q(u_{h}(\bm y)) - Q(P_{M,S}[u_{h,k}](\bm y)) \right|
 \end{equation}
 where $\Xi_{\rm test}$ is $10,000$ uniformly distributed points and $u_{h}$ is some reference finite element solution. In order to generate a reduced basis in our examples we will use a training set $\Xi_{\rm train}$ of $1,000$ uniformly distributed points.
 
{We note for this particular quantity of interest it can be shown using an Aubin-Nitsche duality argument from \cite{CHEN2016470,CHEN201584} that the convergence rates will be twice as great as those in estimate \eqref{conv_rates}.}
 
 \subsection{Example 1}
 
\begin{figure}[h!]
 	\label{fig1}
 	\center
 	\includegraphics[width=8.cm]{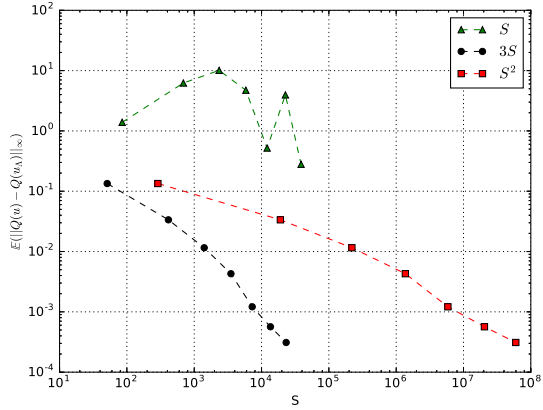} 
 	\caption{Comparing the convergence of the DLS method utilizing different numbers of sample points with respect to the cardinality of the polynomial basis.  }
 \end{figure}

 For our first example, we examine the convergence and stability of the DLS method independent of any reduced basis. Specifically, we are interested in illustrating that a linear rule maintains sufficient stability for the least-squares problem (in a moderately sized dimension) utilizing a quasi-optimal polynomial space. This is motived by computational necessity in larger dimensions for which the cardinality of the quasi-optimal polynomial space can grow very quickly. We choose the  number of sample points {$S = M$}, $3M$, and $M^{2}$. 
 As can be seen in Figure \ref{fig1}, we obtain similar levels of accuracy using the linear rule as we do when using a quadratic rule in agreement with the numerical findings in \cite{CCMFT2013, MFST2013}. {We also see that taking $S = M$ leads to an unstable approximation indicating that some scaling constant is required.}   
 
\subsection{Example 2}

Next, we are interested in the offline and online computational cost   of the DLS method compared to that of the RB-DLS method. 
Beginning with the offline complexity of the DLS method we see, that the majority of the cost is incurred from setting up the right-hand side $\mathbf{U}$ and then solving the system (\ref{lspde}). To solve (\ref{lspde}) we can use any number of methods, two of the more popular being the LU factorization or QR factorization, both of which have the same order of computational complexity. It then follows that the complexity for the DLS algorithm will scale as 
\begin{equation}
\label{costFull}
{\rm DLS}_{\text{cost}} = S \times \mathcal{O}\Big({J^{\alpha}}\Big) + \mathcal{O}(M^{3}) + \mathcal{O}(M^{2}) \times J, 
\end{equation}
where $\mathcal{O}\left({J^{\alpha}}\right)$ is the cost for solving the finite element system where $\alpha$ depends on both the solver and spatial dimension, $\mathcal{O}(M^{3})$ is the cost associated with the LU or QR decomposition, and $\mathcal{O}(M^{2}) \times {J}$ is the cost for solving the system (\ref{lspde}). 

Next, we  analyze the algorithm with the reduced basis incorporated into it. The construction of the reduced basis scales as
\begin{equation}
{\rm RB}_{\text{cost}} = \mathcal{O}(S_{\text{train}}) \times \Big(\sum_{\ell=1}^{\kkk-1} w_{\text{online}}(\ell)\Big) + \kkk \times \mathcal{O}\Big({J^{\alpha}}\Big),
\end{equation}
where $\mathcal{O}(S_{\text{train}})$ is the cost of a max search in our training set, and $w_{\text{online}}(\ell) = \mathcal{O}(\ell^{3})$ is the cost for calculating $\hat{e}_{h,\kkk}$ and $u_{h,k}(\bm y)$ for a value $\bm y \in \Xi_{\rm train}$. The total cost for our algorithm, assuming no online enrichment of the reduced basis is necessary, will thus scale as
\begin{equation}
\label{costReduced}
\text{RB-DLS}_{\text{cost}} = RB_{\text{cost}} + S \times \mathcal{O}({K}^{3}) + \mathcal{O}(M^{3}) + \mathcal{O}(M^{2}) \times {K} + {S \times N^{2} \times K^{2}},
\end{equation}
{where $S \times \mathcal{O}({K}^{3})$ is the cost for solving the reduced basis system for $\{{\bm y}_i\}_{i=1}^S$ to form $\mathbf{W}$ in \eqref{lspde1}, $\mathcal{O}(M^{3})$ is the cost associated with the LU or QR decomposition, $\mathcal{O}(M^{2}) \times K$ is the cost for solving the system \eqref{lspde1} and $S \times N^{2} \times K^{2}$ is the cost for evaluating the error bound $\hat{e}_{h,\kkk}$.   }

We see that the complexity of RB-DLS is dominated by the term $\mathcal{O}(M^{2}) \times K$ when $K$ is large. On the other hand, when ${K}$ is {small} the complexity of both algorithms is dominated by the term $\mathcal{O}(M^{3})$. The key to the computational savings witnessed in the reduced-basis method is that the cost of the reduced-basis algorithm will be independent of ${J}$ except in the offline portion.
As seen in the above discussion, for large values of $J$ the computational cost of the algorithm is dominated by the cost of finite element solves; therefore, we measure the offline computational cost  in terms of the total number of full finite element solves necessary for the construction of the DLS and RB-DLS approximations. 

\begin{figure}[h!]
	\centering
	\includegraphics[width=8.cm]{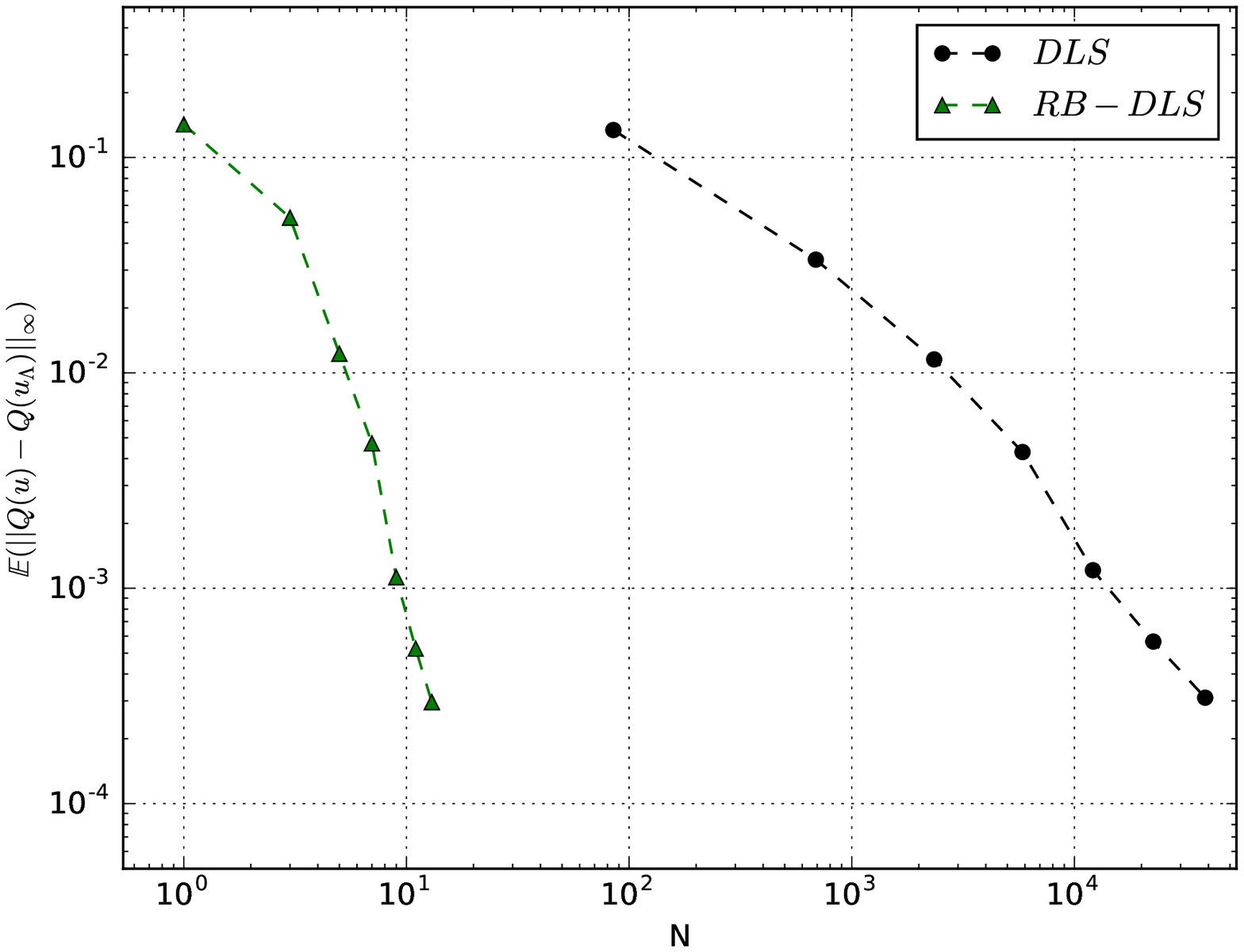}
	\includegraphics[width=8.cm]{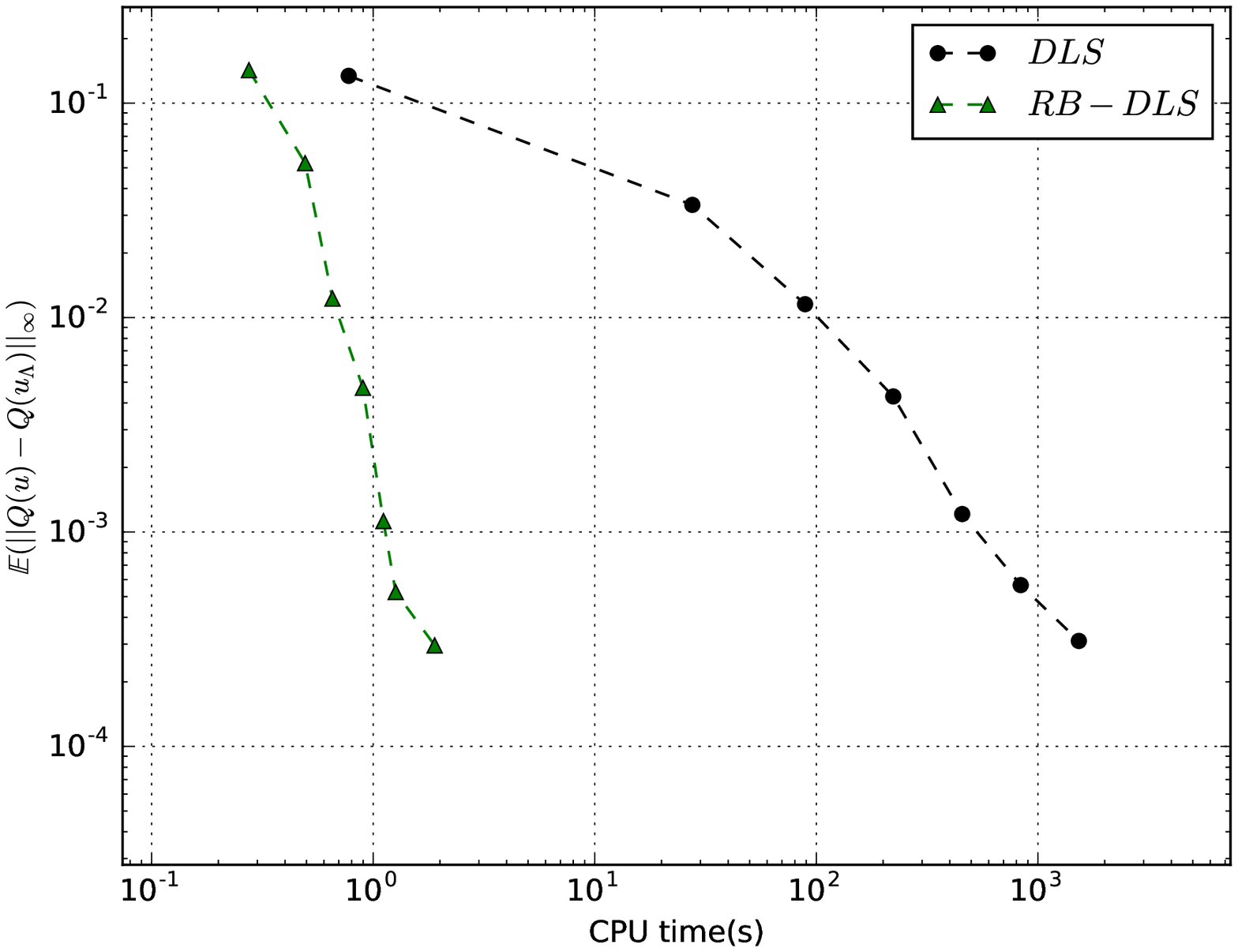}  
	\vspace{-.1in}
	\caption{A comparison of the offline (top) and online (bottom) error versus the cost for the DLS method versus the RB-DLS method {measured in terms of the number of finite element solves in the offline phase (top) and CPU time in the online phase (bottom)}.}
	\label{fig2}
\end{figure}

Turning to the online computational cost, as described in Section \ref{ReducedBasis}, the cost of evaluating DLS for a given $\bm y \in \Gamma$ is of $\mathcal{O}(MJ)$ versus $\mathcal{O}(k(M+J))$ for RB-DLS. To illustrate the significant cost savings of RB-DLS, we compare the total CPU time (in seconds) it takes to compute all RB-DLS and DLS approximations for all $\bm y \in \Xi_{\rm train}$. As shown in Figure \ref{fig2}, we observe significant offline and online computational cost savings while still being able to achieve similar levels of accuracy from the RB-DLS method.

\section{Conclusions}

We integrated a reduced-basis method into the discrete least-squares framework utilizing a new quasi-optimal polynomial space. Through our numerical results, we demonstrated significant cost savings in both the offline and online portions of the discrete least-squares-reduce basis method compared to that for the original discrete least-squares algorithm. Again, we would like to emphasize that reduced basis plays a critical role in solving large-scale UQ problems involving expensive finite element discretization (e.g., with a very fine mesh), especially in the online phase.  We note that this method is not without drawbacks. For the case where the Kolmogorov width of the PDE solution does not decay quickly we will need to use a large number of reduced basis functions in order to obtain an accurate reduced basis approximation. This could potentially make the DLS-RB method more expensive than the standalone DLS method. Additionally the quasi-optimal polynomial basis used in this work only applies to the parametrized diffusion equation. Thus for more complicated PDEs a different polynomial basis would have to be used. This new basis may have significantly worse stability and convergence properties when coupled with the DLS method possible rendering it ineffective. 

{We note that while we paired the reduced basis method with the discrete least-squares algorithm, it is also possible to combine it with sparse-grid method as done in \cite{CHEN2016470,CHEN201584,Chen2014}. A detailed comparison of these approaches has not been done and will be a subject of future research.} 

As shown in this work, a polynomial approximation of the solution map $(\mathbf x, \bm y) \rightarrow u(\mathbf x, \bm y)$ without using a reduced basis may lead to an unaffordable online cost in terms of storage requirement. Thus, model reduction should become a standard procedure in approximating/recovering the solution map $(\mathbf x, \bm y) \rightarrow u(\mathbf x, \bm y)$. 


\bibliographystyle{spmpsci}
\bibliography{myrefs}


\end{document}